\documentclass{amsart}
\usepackage{graphics}

\newcommand{\oct}{\mathcal{O}^{\cdot,D,\cdot}_{\kappa,\cdot,v}(n)}

\newcommand{\cc}{(U,\widetilde{U}/\Gamma, \pi)}
\newcommand{\ccp}{(U',\widetilde{U}'/\Gamma', \pi')}

\newcommand{\wtu}{\widetilde{U}}

\newcommand{\al}{\alpha}
\newcommand{\kp}{\kappa}
\newcommand{\ga}{\gamma}
\newcommand{\Ga}{\Gamma}
\newcommand{\ld}{\lambda}
\newcommand{\om}{\omega}
\newcommand{\wtom}{\widetilde{\omega}}
\newcommand{\pt}{\frac{\pi}{2}}
\newcommand{\apq}{\mathfrak{d}_{pq}}
\newcommand{\aqp}{\mathfrak{d}_{qp}}

\newcommand{\apK}{\mathfrak{d}_{\scriptscriptstyle{pK}}}
\newcommand{\fa}{\mathfrak{a}}
\newcommand{\fb}{\mathfrak{b}}
\newcommand{\fc}{\mathfrak{c}}
\newcommand{\tfa}{\tilde{\mathfrak{a}}}
\newcommand{\tfb}{\tilde{\mathfrak{b}}}

\newcommand{\tv}{\tilde{v}}
\newcommand{\tw}{\tilde{w}}
\newcommand{\tp}{\tilde{p}}
\newcommand{\tq}{\tilde{q}}
\newcommand{\tif}{\tilde{f}}

\newcommand{\op}{\overline{p}}
\newcommand{\oq}{\overline{q}}
\newcommand{\ve}{\varepsilon}
\newcommand{\vp}{\varphi}
\newcommand\vol{\operatorname{Vol}}

\newcommand{\ch}{\mathcal{H}}

\newcommand{\cs}{\mathcal{S}}
\newcommand{\ca}{\mathcal{A}}
\newcommand\e{\operatorname{exp}}
\newcommand\seg{\operatorname{seg}}
\newcommand\snk{\operatorname{sn_{\kappa}}}

\newcommand{\scc}{M^n_{\kappa}}
\newcommand{\ccb}{B^n_{\kp}}
\newcommand{\lh}{\infty}
\newcommand{\parttxi}{\frac{\partial}{\partial \tilde{x}^i}}

\newcommand{\sings}{O - \Sigma_O}

\newcommand\inr{\operatorname{int}}
\newcommand\Rc{\operatorname{Ric}}

\newtheorem{theorem}{Theorem}[section]

\newtheorem{lemma}[theorem]{Lemma}
\newtheorem{prop}[theorem]{Proposition}
\newtheorem{defn}[theorem]{Definition}

\newcommand{\lemref}[1]{Lemma~\ref{#1}}
\newcommand{\propref}[1]{Prop\-o\-si\-tion~\ref{#1}}

\newcommand{\figref}[1]{Fig\-ure~\ref{#1}}

\begin{document}
\title[Spectral bounds on orbifold isotropy]{Spectral bounds on orbifold 
isotropy}
\author[E. Stanhope]{Elizabeth Stanhope}
\address{Willamette University, Mathematics Department, 900 State Street, Salem OR 97301}


\thanks{{\it Math. Classification:} 58J50 \ 53C20}

\email{estanhop@willamette.edu}

\bibliographystyle{alpha} 

\maketitle

\section*{Introduction}
    
An underlying theme in differential geometry is to uncover
information about the topology of a Riemannian manifold
using its geometric structure.  The present investigation carries this
theme to the category of Riemannian orbifolds.  In particular we ask:
If a collection of isospectral orbifolds satisfies a uniform
lower bound on Ricci curvature, do orbifolds in the collection
have similar topological features?  Can we say more if we require the
collection to satisfy a uniform lower sectional curvature bound?  We will 
assume throughout that all orbifolds are connected and closed.

Our inquiry begins with a review of the fundamentals of doing geometry
on orbifolds in Sections \ref{orbifold background} through \ref{geod}.
Riemannian orbifolds, first defined by Satake in \cite{MR18:144a}, are spaces 
that are locally modelled
on quotients of Riemannian manifolds by finite groups of isometries.
These sections examine topics including the behavior of geodesics on orbifolds, 
and integration on
orbifolds.

In Section \ref{ofdcompgeom} we see how the geometry of 
orbifolds
with lower curvature bounds can be studied by comparing them to manifolds
with constant curvature.  This section builds on the work in
\cite{MR94d:53053}.

The eigenvalue spectrum of the Laplace operator on an orbifold is
  introduced in Section \ref{spec intro}.  We confirm that several familiar 
spectral theory tools from the manifold
  setting carry over to orbifolds.  An orbifold version of Weil's asymptotic
formula from \cite{MR2001k:58060} is stated, showing that the dimension and 
volume of an
  orbifold can be deduced from its spectrum.

The last two sections develop the proofs of two
affirmative answers to our main questions.  In both statements below
we assume the orbifolds under consideration are compact and orientable.

\smallskip

\noindent \bf{Main Theorem 1:} \it   Let $\cs$ be a collection of isospectral 
Riemannian orbifolds that share a uniform lower bound
$\kp(n-1)$ on Ricci curvature, where $\kp \in \mathbf{R}$.  Then there are only 
finitely
many possible isotropy types, up to isomorphism, for points in an
orbifold in $\cs$. \rm

\smallskip

\noindent \bf{Main Theorem 2:} \it Let $isol \cs$ be a collection of 
isospectral Riemannian orbifolds with only
isolated singularities, that share a uniform lower bound $\kp \in
\mathbf{R}$ on sectional curvature.  Then there is an upper
bound on the number of singular points in any orbifold, $O$, in $isol  \cs$
depending only on $Spec(O)$ and $\kp$. \rm

\smallskip

\noindent Note that there exist examples of constant curvature one isospectral orbifolds 
which possess distinct isotropy.  Thus Main Theorem 1 cannot be improved to 
uniqueness.

The proofs of these results break down into two steps.  The
first step is to convert spectral information into explicit bounds
on geometry.  As mentioned above, the dimension and the volume of an
orbifold are determined by its spectrum.  In Section \ref{chengchapt}
we obtain an upper bound on the diameter of an orbifold which depends
only on the orbifold's spectrum, and the presence of a lower bound on
Ricci curvature.  The technique used to derive this
diameter bound parallels a similar one from the manifold setting given
in \cite{MR93f:53034}.  The main ingredient used is an orbifold version of
Cheng's Theorem.  

The second step in proving these theorems is to examine families of $n$-orbifolds that satisfy an
upper diameter bound, and lower bounds on curvature and
volume.  By the work in the first step, results that hold for these
families also hold for families of isospectral orbifolds having a
uniform lower bound on curvature.  The first main theorem is shown
using volume comparison techniques.  The second main theorem relies both on
tools from comparison geometry, and on a careful analysis of the orbifold
distance function, generalizing results of Grove and Petersen
\cite{MR90a:53044} to the orbifold setting.  This analysis is the focus of 
Section \ref{main lemma chap}.

\smallskip

\noindent \bf{Acknowledgements.} \rm
The author would like to thank her thesis advisor, Carolyn S. Gordon,
for her guidance and patience during the course of this work.

\section{Smooth Orbifolds} \label{orbifold background}

An orbifold is a generalized manifold
arrived at by loosening the requirement that the space be locally modelled on 
$\mathbf{R}^n$, and instead
requiring it to be locally modelled on $\mathbf{R}^n$ modulo the action
of a finite group.  This natural generalization allows orbifolds to
possess `well-behaved' singular points.  In this section we make these
ideas precise and set up some basic tools that will be used throughout this 
text.

We first recall the definition of smooth orbifolds given by Satake in
\cite{MR18:144a} and \cite{MR20:2022}.  In order to state the definition we 
need to
specify what is meant by a chart on an orbifold, and what it means to
have an injection between charts.

\begin{defn}  Let $X$ be a Hausdorff space and $U$ be an open set in
  $X$.  An orbifold
    coordinate chart over $U$, also known as a uniformizing
  system of $U$, is a triple $\cc$ such that:  
\begin{enumerate}
\item $\widetilde{U}$ is a connected open subset of $\mathbf{R}^n$,
\item  $\Gamma$ is a finite group of diffeomorphisms acting on
  $\widetilde{U}$ with fixed point set of codimension $\ge$ 2, and 
\item $\pi:\widetilde{U} \rightarrow U$ is a
  continuous map which induces a homeomorphism between $\wtu/\Ga$ and $U$,
  for which $\pi \circ \gamma = \pi$ for all $\gamma \in \Gamma$. 
\end{enumerate}
\end{defn}

Now suppose $X$ is a Hausdorff space containing open subsets $U$ and
$U'$ such that $U$ is contained in $U'$.  Let $\cc$ and $\ccp$ be charts over 
$U$ and $U'$, respectively.
\begin{defn}  An injection $\ld: \cc \hookrightarrow \ccp$ consists of an
  open embedding $\ld:\wtu \hookrightarrow \wtu '$ such that $\pi =
  \pi' \circ \ld$, and for any $\gamma \in
  \Gamma$ there exists $\gamma' \in \Gamma'$ for which $\ld \circ
  \gamma = \gamma' \circ \ld$.
\end{defn}
Note that the correspondence $\gamma \mapsto \gamma'$ given above defines an
injective homomorphism of groups from $\Gamma$ into $\Gamma'$.

\begin{defn}\label{orbifold}  A smooth orbifold $(X, \mathcal{A})$ consists of 
a Hausdorff space
  $X$ together with an atlas of charts $\mathcal{A}$ satisfying the following 
conditions:
\begin{enumerate}
\item For any pair of charts $\cc$ and $\ccp$ in $\mathcal{A}$ with
  $U \subset U'$ there exists an injection $\ld: \cc \hookrightarrow \ccp$. 
\item The open sets $U \subset X$ for which there exists a chart $\cc$
  in $\mathcal{A}$ form a basis of open sets in $X$.
\end{enumerate}
\end{defn}
Given an orbifold $(X,\ca)$, the space $X$ is referred to as the
\emph{underlying space} of the orbifold.  Henceforth specific reference to an 
orbifold's underlying space and atlas
of charts will be dropped and an orbifold $(X,\ca)$ will be
denoted simply by $O$.  

Take a point $p$ in an orbifold $O$ and let $\cc$ be a coordinate chart
about $p$.  Let $\tp$ be a point in $ \wtu$ such that $\pi(\tp) = p$ and let 
$\Gamma_{\tp}^U$ denote the isotropy
group of
$\tp$ under the action of $\Gamma$.  It can be shown that the group
$\Gamma_{\tp}^U$ is actually independent of both
the choice of lift and the choice of chart (see \cite{MR94d:53053}), and so can 
sensibly
be denoted by $\Gamma_p$.  We call $\Gamma_p$ the \emph{isotropy group
of $p$}.  Points in $O$ that have a non-trivial isotropy group
are called \emph{singular points}.  We will let $\Sigma_O$ denote the set
of all singular points in $O$.  

Before describing more properties of orbifolds, we state a proposition
which gives an important class of
orbifolds. A proof can be found in \cite{thur}.
\begin{prop}
Suppose a group $\Gamma$ acts properly discontinuously on a manifold
$M$ with fixed
point set of codimension greater than or equal to two.  Then the quotient space 
$M/\Ga$ is an orbifold.
\end{prop}

An orbifold is called \emph{good}
(\emph{global} is also used) if it
arises as the quotient of a manifold by a properly discontinuous group
action.  Otherwise the orbifold is called \emph{bad}.

Suppose $O = M/\Ga$ is a good orbifold.  We
can extend the action of $\Ga$ on $M$ to an action on $TM$ by setting
$\gamma . (\tp,v) = (\ga (\tp), \ga_{*_{\tp}} v)$ for all $\ga \in \Ga$ and 
$(\tp,v)
\in TM$.  The quotient of $TM$ by this new action is the tangent
bundle, $TO$, of the orbifold $O$.  For $\tp \in M$ let $p \in O$ be the
image of $\tp$ under the quotient.  By taking the differentials at
$\tp$ of
elements of the isotropy group of $p$, we form a
new group that acts on $T_{\tp}M$.  Because this group is independent
of choice of lift, we can denote it by
$\Ga_{p*}$.  The fiber in $TO$
over $p$ is $T_{\tp}M/\Ga_{p*}$, and is denoted $T_pO$.  Because $T_pO$ need 
not be a vector space, it is called the \emph{tangent cone to $O$ at $p$}.

Locally all orbifolds are good, so the construction above gives a
local way to work with tangent cones to orbifolds.  A full
construction of orbifold tangent bundles, as well as general bundles over
orbifolds, is given in \cite{MR20:2022}. 

\section{Riemannian Metrics on Orbifolds}\label{Riemofds}

After giving the definition of smooth functions on orbifolds, we move
on to more general tensor fields including the Riemannian metric.  In this section, 
and all that follow, we
will assume that each orbifold has a second countable underlying
space.  In addition to \cite{MR18:144a} and \cite{MR20:2022}, useful references 
for this material include \cite{MR94d:53053} and \cite{MR94c:58040}. 

\begin{defn}A map $f:O \rightarrow \mathbf{R}$ is a smooth function on $O$ if 
on each
chart $\cc$ the lifted function $\tif = f \circ \pi$ is a smooth
function on $\wtu$.  
\end{defn}

\begin{defn}  Let $\cc$ be an orbifold coordinate chart.  
\begin{enumerate}
\item For any tensor field $\widetilde{\omega}$ on $\wtu$
precomposing by $\ga \in \Ga$ gives a new tensor field on $\wtu$, denoted
$\widetilde{\omega}^{\ga}$.  By averaging in this manner we obtain a
$\Ga$-invariant tensor field, denoted $\widetilde{\omega}^{\Ga}$, on $\wtu$:
\begin{equation*}
\widetilde{\omega}^{\Ga} = \frac{1}{|\Ga|}\sum_{\ga \in \Ga} \widetilde{\omega}^
{\ga}.
\end{equation*}
Such a $\Ga$-invariant tensor field on $\wtu$
gives a tensor field $\omega$ on $U$.  
\item A smooth tensor field on an orbifold is one
that lifts to smooth tensor fields of the same type in all local covers.
\end{enumerate}
\end{defn}

A Riemannian metric is obtained on a good orbifold, $O = M/\Ga$, by
specifying a Riemannian metric on $M$ that is invariant under the
action of $\Ga$.  This also gives a local notion of Riemannian metric
which leads to the definition of a Riemannian metric for general orbifolds.
Let $O$ be a general orbifold and $\cc$ be one of its coordinate charts.  
Specify a Riemannian metric $g^{\wtu}$ on $\wtu$.  By averaging as above we can 
assume that
this metric is invariant under the local group action, and so gives
 a Riemannian metric $g^U$ on $U$.  Now do this for each chart of
 $O$.  By patching the local metrics together using a partition of
 unity, we obtain a global Riemannian metric $g$ on $O$.  A smooth
 orbifold together with a Riemannian metric is called a \emph{Riemannian 
orbifold}.

In the construction above, the Riemannian metric $g^{\wtu}$ on $\wtu$ is 
invariant under the action of $\Gamma$.  Another way
to say this is that locally Riemannian orbifolds look like the quotient
of a Riemannian manifold by a finite group of isometries.  By a suitable choice 
of coordinate charts (see \cite{MR94c:58040},
p. 318) it can be assumed that the local group actions are by finite
subgroups of $O(n)$ for general Riemannian orbifolds, and finite subgroups of 
$SO(n)$ for
orientable Riemannian orbifolds.

Objects familiar from the Riemannian geometry of manifolds are defined
for orbifolds by using the Riemannian metrics on the local covers.  For 
example, we say
that a Riemannian orbifold $O$ has sectional curvature
bounded below by $k$ if every point is locally covered by a
manifold with sectional curvature greater than or equal to $k$.  Ricci
curvature bounds are defined similarly.  We define angles in the following 
manner.

\begin{defn}  Let $p$ be a point in a Riemannian orbifold that lies in a
  coordinate chart $\cc$.  Take $\tp$ to be a lift of $p$ in
  $\wtu$.  For vectors $v$ and $w$ in $T_pO$ let $\tv_1, \tv_2, \dots, \tv_r$
  denote the set of lifts of $v$, and $\tw_1, \tw_2, \dots, \tw_s$
  denote the set of lifts of $w$, in $T_{\tp}\wtu$.  The angle between $v$ and 
$w$ in $T_pO$ is defined to be, 
\begin{align*}
\angle(v,w) = \min_{\substack{ 
                       i = 1, 2 \dots, r \\
                       j = 1, 2 \dots, s 
                   }}   \{\angle(\tv_i, \tw_j)\}.
\end{align*}
\end{defn}

If $O = M/\Ga$ is a good Riemannian orbifold, the quotient of the unit
tangent bundle of $M$ by $\Gamma$ yields the unit tangent bundle of
the orbifold, $SO$.  The \emph{unit tangent cone to $O$ at $p$},
denoted $S_pO$, is the fiber over $p$ in this bundle.  Alternatively
the unit tangent cone is the set of all unit vectors in $T_pO$.

A particularly useful type of chart about a point $p$ in a Riemannian orbifold
is one for which the group action is by the isotropy group
of $p$.  This type of chart is called a \emph{fundamental coordinate
chart} about $p$.  Every point in a Riemannian orbifold lies in a
fundamental coordinate chart  (see \cite{MR94d:53053}, p. 40).

\section{Geodesics and Segment Domains for Orbifolds}\label{geod}

We now examine the structure of geodesics in orbifolds.  In this
discussion, length minimizing geodesics will be referred to as segments.

Let $p$ be a point in a Riemannian orbifold $O$, and let $\cc$ be a
coordinate chart about $p$.  For every $v \in
S_pO$ there is a segment $\gamma_v$ that emanates from $p$ in the
direction of $v$.  To see this, take $\tp$ to be a lift of $p$ in
$\wtu$, and $\tv$ to be a lift of
$v$ in $S_{\tp}\wtu$.  For small $t$ we have the
segment $\e_{\tp}t\tv$ emanating from $\tp$ in $\wtu$.  The image of
this segment under $\pi$ is a segment in $O$ that leaves $p$ in the direction
of $v$.  Thus within a coordinate chart about $p$ we can define the
exponential map, $\e_p tv$, by projecting $\e_{\tp}t\tv$ to $U$.  Note
that this definition is well-defined as it is independent of choice of lift.

To obtain the exponential map globally on an orbifold we extend these
locally defined geodesics as far
as possible.  More precisely, for $v \in S_pO$ let $\gamma_v(t)$
denote the geodesic emanating from $p$ in the direction $v$.  Then for
all $t_0 \in [0,+\lh)$ where $\ga_v(t_0)$ is defined, set $\e_p t_0v =
\ga_v(t_0)$.  

In Proposition 15 of \cite{MR94d:53053} it is shown that if a segment is not
entirely contained within the singular set, it can only intersect the
singular set at its end points.  So a segment that contains any
manifold points must stop when it hits the singular set.
Consequentially if an orbifold is to be geodesically
complete, no obstruction by singular points can occur.  Thus the
singular set of a geodesically complete orbifold must be empty,
implying the orbifold is actually a manifold.  In what follows the
word complete will be used to describe orbifolds that, together with their 
distance
functions, are complete as metric spaces.  An analogue of the
Hopf-Rinow Theorem for length spaces (see \cite{MR2000d:53065}, p. 9) implies 
that if an orbifold is
complete, then any two points in the orbifold can be joined by a
segment.  

Suppose $O$ is a complete orbifold and consider the manifold obtained
by excising its singular set, $O - \Sigma_O$.  The preceding
observations imply that any two points in $O - \Sigma_O$ are connected
by a segment that lies entirely within $O - \Sigma_O$.  Thus we see
that $O - \Sigma_O$ is a convex manifold.  This fact will be used
extensively in what follows.

We will now consider the segment domain of an orbifold.  

\begin{defn}\label{segdom}
The segment domain of a point $p$ in an orbifold $O$ is denoted
by $\seg (p)$ and is defined as follows:
\begin{equation*}
\seg (p) = \{ v \in T_pO : \e_p tv : [0, 1] \rightarrow O \text{ is a
  segment} \}.
\end{equation*}
The interior of the segment domain of $p$, $\seg^0(p)$, is defined by:
\begin{equation*}
\seg^0(p) = \{ vt : t \in [0, 1), v \in \seg(p) \}.
\end{equation*}
\end{defn}

For $p \in O$, the image of the boundary of $\seg(p)$ under the exponential
map at $p$ is called the \emph{cut locus} of $p$ in $O$.  The cut locus of
$p$ is denoted by $\text{cut}(p)$.  This set consists of the points in
$O$ beyond which geodesics from $p$ first fail to minimize distance.

The use of the segment domain in what follows relies on the following
lemma.  Its proof is analogous to that of the manifold case.

\begin{lemma}\label{expadfo}Let $O$ be a complete Riemannian orbifold and take 
$p \in
  \sings$.  Then $\exp_p: \seg^0(p) \rightarrow O$ is a diffeomorphism
  onto its image.
\end{lemma}

We end this section by defining integration on orbifolds
and by describing a useful integration technique.  Suppose that $O$ is a compact orientable Riemannian orbifold.  Let
$\omega$ be an $n$-form on $O$ such that the support of $\omega$ is contained 
in the chart $\cc$.  We define the integral of
$\omega$ over $O$ as follows,
\begin{equation*}
\int_O \om = \frac{1}{|\Ga|} \int_{\wtu} \wtom , 
\end{equation*}
where $\wtom = \om \circ \pi$.
By using the injections provided by the orbifold structure, one can
check that this definition does not depend on the choice of coordinate
chart.  The integral of a general $n$-form is defined using a partition of 
unity,
as in the manifold case.

Sometimes it will be more convenient to compute integrals using the
  following technique.  Let $p \in \sings$.  Then $p$ has a manifold 
neighborhood in $O$ upon
  which we can consider the usual manifold polar coordinates.  The
  volume density in these polar coordinates is
  $\sqrt{\text{det}(g_{\alpha \beta}(r, \theta))}$, which will be denoted
  by $\rho(r,\theta)$ for convenience.

\begin{prop}\label{intlemma}  Let $O$ be a complete Riemannian orbifold, with 
$p \in
  \sings$ and suppose $f \in C^{\lh}(O)$.  Then,
\begin{equation*}
\int_O f\ dV = \int_{\seg^0(p)} (f \circ \exp_p) \rho(r, \theta) \
dr d\theta.
\end{equation*}
\end{prop}

\section{Comparison Geometry Background}\label{ofdcompgeom}
The geometry of hyperbolic space, Euclidean space and the sphere is
very well developed, in contrast to that of manifolds with variable
curvature.  The idea behind comparison geometry is to study spaces
with variable curvature by comparing them to the simply connected spaces
with constant sectional curvature.  

In this section we confirm that several familiar comparison results are valid 
in the orbifold setting.  The following notation will be helpful.  We will use 
$\scc$ to denote the
simply connected $n$-dimensional space form of constant curvature
$\kp$.  The open $r$-ball in $\scc$ will be denoted by
$B^n_{\kp}(r)$.  As in Section \ref{geod}, the volume
density of a manifold will be written in polar coordinates as $\rho(r,
\theta)$.  We denote the volume
density on $\scc$ by $(\snk(r))^{(n-1)}$, where $\snk(r)$ is given by:
\begin{equation*}
\snk(r) = \begin{cases}
\frac{\sin{\sqrt{\kp}r}}{\sqrt{\kp}} & \kp > 0 \\
r & \kp = 0 \\
\frac{\sinh{\sqrt{-\kp}r}}{\sqrt{-\kp}} & \kp < 0. 
\end{cases}
\end{equation*}

The Relative Volume Comparison Theorem is generalized to orbifolds in \cite
{MR94d:53053}.

\begin{prop}\label{ofdrvct}(Orbifold Relative Volume
  Comparison Theorem)  Let $O$ be a complete Riemannian orbifold with
  $\Rc(M) \ge (n-1)\kp$.  Take $p \in O$.  Then the function,
\begin{equation*}
r \mapsto \frac{\vol B(p,r)}{\vol \ccb(r)}
\end{equation*}
is non-increasing and has limit equal to $\frac{1}{|\Gamma_p|}$ as $r$ goes to 
zero.
\end{prop}

\noindent Note that this theorem implies a volume comparison theorem for balls in 
orbifolds.  To see this observe that if
  $0 \le r \le R$ then by the theorem above,
\begin{equation*}
\frac{\vol B(p,r)}{\vol \ccb(r)} \ge \frac{\vol B(p,R)}{\vol \ccb(R)}.
\end{equation*}
Taking the limit of this inequality as $r$ goes to zero shows that the
volume of an $R$-ball in $O$ is less than or equal to the volume of an
$R$-ball in $\scc$.

We next specify what is meant by a cone in an orbifold.

\begin{defn}\label{orbifold cone}
Let $p \in O$ and $\fa \subset S_pO$, the tangent sphere to $O$ at
$p$.  The $\fa$-cone at $p$ of radius $r$ is defined to be,
\begin{equation*}
B^{\fa}(p, r) = \{\exp_ptv : (t,v) \in Domain(\e_p), \ t < r, v \in \fa \} 
\subset O.
\end{equation*}  
The associated cone in $T_pO$ is defined as follows, 
\begin{equation*}
B^{\fa}(0, r) = \{tv : (t,v) \in Domain(\e_p), \ t < r, v \in \fa \} \subset 
T_pO.
\end{equation*}
\end{defn}

\begin{figure}[b]
\scalebox{.5}{\includegraphics{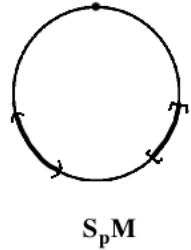}}
\caption{Subset of Unit Tangent Circle}\label{ainspm}
\end{figure}

We illustrate this definition in the case of surfaces.  In \figref
{ainspm} a subset of
the unit tangent circle at a point $p$ in a surface $M$ is specified.
\begin{figure}
\scalebox{.5}{\includegraphics{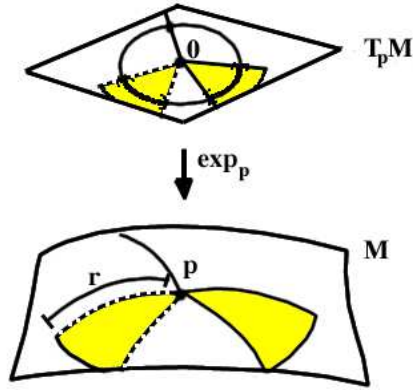}}
\caption{Illustration of Cones}\label{coneil}
\end{figure}
The associated cones of radius $r$ in the tangent space and in the surface are
illustrated in Figure \ref{coneil}.

In Chapter 9 of \cite{MR98m:53001} a volume comparison theorem for cones in 
manifolds is considered.  We will need a version of this theorem that is valid 
for orbifolds.  In order to state this theorem, we will use
the following notation.  We suppose $p$ is
a point in an orbifold $O$ with fundamental coordinate chart $\cc$.  For $A 
\subset
T_pO$, the set $\{\tv \in T_{\tp}\wtu : \pi_{*\tp} \tv \in A\}$ is denoted by 
$\widetilde{A}$.

\begin{prop}\label{ofdcvct}{(Volume comparison theorem for cones in orbifolds.)}
Let $O$ be a complete Riemannian orbifold with $\Rc(O) \ge (n-1)\kp$,
  and take $p \in O$.
  If $\kp > 0$ suppose $r \le
  \pi/\sqrt{\kp}$, otherwise let $r$ be any non-negative real.
  Suppose $\fa$ is an open subset of $S_pO$ with boundary of measure zero,
  and $\op \in \scc$.  Let $I$ be an isometry from $S_{\tp}\wtu$ to
$S_{\op}\scc$, relative to the canonical metric on the unit sphere.  Then,
\begin{equation*}
\vol B^{\fa}(p,r) \le \frac{1}{|\Ga_p|}\vol B^{I(\tfa)}(\op,r).
\end{equation*}
\end{prop}

\begin{proof}
First suppose that $p$ is a manifold point in $O$.  Using the fact that $O - \Sigma_O$ is a 
convex manifold, and that $p$
has trivial isotropy, we conclude:
\begin{equation*}
\vol B^{\fa}(p,r) \le \vol B^{I(\fa)}(\op,r) = \frac{1}{|\Ga_p|} \vol B^{I
(\tfa)}(\op,r).
\end{equation*}

Now suppose $p$ is a singular point in $O$.  Let $(U,\wtu/\Gamma,\pi)$ be a 
fundamental coordinate chart about $p$.  Suppose $\tp \in
\wtu$ projects to $p$, and lift $\fa$ to $\tfa
\subset S_{\tp}\wtu$.  

Choose a vector $v \in \fa$ that points out of the singular set.  Fix
a lift $\tv$ of $v$ in $\tfa$.  Recall that the Dirichlet fundamental
domain centered at $\tv$ of the action of $\Gamma_{*_{\tp}}$ on
$S_{\tp}\wtu$ is the set $\{u \in S_{\tp}\wtu : d(u,\tv) \le
d(u, \gamma_{*_{\tp}} \tv) \ \text{for all} \ \gamma_{*_{\tp}} \in
\Gamma_{*_{\tp}}\}$.  Let $\tfb$ denote the intersection of this
Dirichlet fundamental domain with $\tfa$.  Let $\gamma_{\tv}:[0,\varepsilon) 
\rightarrow \wtu$ be
a portion of the geodesic emanating from $\tp$ in the direction
$\tv$.  Let $\gamma_v$ be the image of $\gamma_{\tv}$ under $\pi$.
Shrink $\ve$ as needed to ensure that $\ga_v([0,t])$ is minimizing for
all $t \in [0,\ve)$ and that $\ve < r$.

The parallel transport map $P_{0,t}:T_{\tp}\wtu \rightarrow
T_{\gamma_{\tv}(t)}\wtu$ is a vector space isometry.  Let
  $\tfb(t)$ be the subset $P_{0,t}(\tfb) \subset S_{\gamma_{\tv}(t)}\wtu$.  Note
    here that $\tfb(0) = P_{0,0}(\tfb) = \tfb$.  This process smoothly
    spreads $\tfb$ along the spheres tangent to points on the geodesic
    $\gamma_{\tv}(t)$.

Using this, for $t \in (0,\ve)$
we can specify a subset $\fb(t)$ of $S_{\gamma_{v(t)}}O$ by $\fb(t) = \pi_
{*\gamma_{\tv}(t)}(\tfb(t))$.  

For $A \subset O$ let $\chi_A$ denote the characteristic function of
$A$ given by:
\begin{equation*}
\chi_A(x)=\begin{cases}
1& x \in A\\
0& x \in O-A
\end{cases}.
\end{equation*}
We will show that as $t$ goes to zero in $[0, \ve)$, the functions $\chi_{B^{\fb
(t)}(\gamma_v(t),r-t)} \rightarrow
\chi_{B^{\fa}(p,r)}$ pointwise a.e.  To do this
we need to check that points in the cone $B^{\fa}(p,r)$ also lie in
nearby cones $B^{\fb(t)}(\gamma_v(t),r-t)$, and points outside of
$B^{\fa}(p,r)$ also lie outside nearby cones
$B^{\fb(t)}(\gamma_v(t),r-t)$.  Because the property of being in a
particular cone depends on distance and angle, we check each of these
in the two cases.

Fix $x$ in the $r$-ball about $p$.  Then, for this $x$, we can find
a $\delta_1> 0$ sufficiently small so that  $x$ will be in the
balls $B(\ga_v(t), r-t)$ for all $t \in [0,\delta_1)$.

Now consider the directions from points on $\gamma_v$ to $x$.  Let $\sigma_t$ 
denote the geodesic from
$\gamma_v(t)$ to $x$.  The fact that $x$ lies in $B^{\fa}(p,r)$ implies that 
$\sigma_0'(0)
\in \fa$.  Noting that $\fa = \fb$ is an open subset of $S_pO$, we can assume 
there is a small neighborhood $\fc$ about $\sigma_0'(0)$ in
$\fb$.  By lifting and translating as above we have $\fc(t) \subset
\fb(t)$ for $t\in (0,\ve)$.  By continuity, $\sigma_t'(0)$ will remain
in $\fc(t)$ for $t$ small, say for $t \in [0,\delta_2)$.  Thus
$\sigma_t'(0)$ will remain in $\fb(t)$ for $t \in [0, \delta_2)$.

Set $\delta = \min\{\delta_1, \delta_2\}$.  The previous two paragraphs
imply that $x$ lies in the cones $B^{\fb(t)}(\ga_v(t),r-t)$ for $t \in [0, 
\delta)$.

Now suppose that $x$ lies outside of the cone $B^{\fa}(p,r)$.  This
means that either the distance between $p$ and $x$ is larger than $r$,
or the direction from $p$ to $x$ lies outside of $\fa$.  We
need to confirm that in either of these cases, $x$ also lies outside
of cones $B^{\fb(t)}(\ga_v(t),r-t)$ for small $t$.  Because the balls $B(\ga_v
(t), r-t)$ are contained in $B(p,r)$, it is
clear that points lying outside of $B(p,r)$ are also outside of
$B(\ga_v(t), r-t)$ for $t \in [0,\ve)$.

Suppose that $x$ fails to be in $B^{\fa}(p,r)$ because the
direction from $p$ to $x$ is not in $\fa$.  As before let $\sigma_t$ denote the 
geodesic from
$\ga_v(t)$ to $x$.  That the direction from $p$ to $x$ lies
outside of $\fa$ is written more precisely as $\sigma_0'(0) \in S_pO -
\fa$.  Disregarding points on the boundary of $B^{\fa}(p,r)$, we
can assume the slightly stronger statement that $\sigma_0'(0) \in S_pO -
\overline{\fa}$.  Take a small neighborhood $\fc$ about $\sigma_0'(0)$
in $S_pO - \overline{\fa}$.  By continuity there is an $\eta >0$
such that $\sigma_t'(0)$ lies outside of $\fb(t)$ for all $t \in
[0,\eta)$.  Thus $x$ lies outside of the cones
$B^{\fb(t)}(\ga_v(t),r-t)$ for $t\in [0,\eta)$ as desired.

We are now able to apply the Lebesgue
Dominated Convergence Theorem to obtain,
\begin{equation*}
\vol B^{\fb(t)}(\gamma_v(t),r-t) = \int_{O - \Sigma_O}
\chi_{B^{\fb(t)}(\gamma_v(t),r-t)} dV \stackrel{t \rightarrow 0}
{\longrightarrow} \int_{O - \Sigma_O}
\chi_{B^{\fa}(p,r)} dV = \vol B^{\fa}(p,r).
\end{equation*}

For $t \in (0, \ve)$, each $\gamma_v(t)$ is a manifold point in $O$.
Because the proposition holds for manifold points, we conclude that if $t \in 
(0,\ve)$ then,
\begin{equation*}
\vol B^{\fb(t)}(\gamma_v(t),r) \le \vol B^{I(\fb(t))}(\op,r).
\end{equation*}  
Now on $(0,\ve)$ we have $\fb(t)$ isometric to $\tfb(t)$ via  $\pi_{*\gamma_
{\tv}(t)}$, and $\tfb(t)$
is isometric to $\tfb$ via $P_{0,t}$.  Thus,
\begin{equation*}
\vol B^{\fb(t)}(\gamma_v(t),r) \le \vol B^{I(\tfb)}(\op,r).
\end{equation*}  

Taking the limit as $t \rightarrow 0$ in this inequality yields,
\begin{equation*}
\vol B^{\fa}(p,r) \le \vol B^{I(\tfb)}(\op,r).
\end{equation*}  
Finally because the translates of $\tfb$ cover $\tfa$ and overlap on a
set of measure zero, we have,
\begin{equation*}
\vol B^{I(\tfb)}(\op,r) = \frac{1}{|\Ga_p|} \vol B^{I(\tfa)}(\op,r).
\end{equation*}
\end{proof}

We end this section with a version of Toponogov's Theorem for orbifolds.  In 
\cite{MR94d:53053} it is shown
that orbifolds with sectional curvature bounded below by $\kp \in
\mathbf{R}$ have Toponogov curvature greater than or equal to $\kp$ in
the sense of length spaces.  In particular, an orbifold with a lower
bound $\kp$ on sectional curvature is an Alexandrov space with
curvature bounded below by $\kp$.

The following proposition is proven in \cite{MR96c:53064}.

\begin{figure}[t]
\scalebox{.45}{\includegraphics{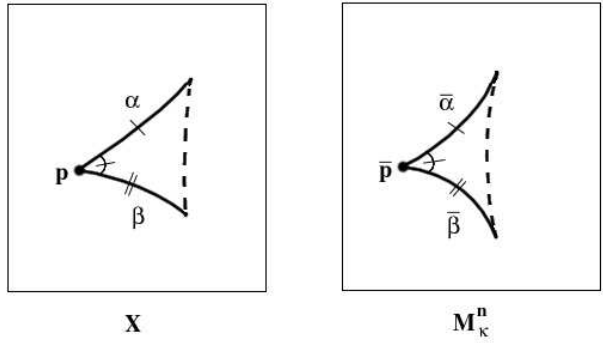}}
\caption{Hinges in Toponogov's Theorem}\label{top}
\end{figure}

\begin{prop}  Let $X$ be an Alexandrov space with curvature bounded
  below by $\kp$.  Let $\al:[0,a] \rightarrow X$ and $\beta:[0,b]
  \rightarrow X$ be geodesics with $\al(0) = \beta(0) = p$  (see
  \figref{top}).  Let
  $\overline{\al}$ and $\overline{\beta}$ be geodesics from point
  $\op$ in $\scc$ with the same lengths as $\al$ and $\beta$,
  respectively, and with $\angle(\al'(0),\beta'(0)) =
  \angle(\overline{\al}'(0),\overline{\beta}'(0))$.  Then
  $d(\al(a),\beta(b)) \le d(\overline{\al}(a),\overline{\beta}(b))$.
\end{prop}

We conclude that the hinge version of Toponogov's Theorem is valid for
orbifolds.

\section{Spectral Geometry Background}\label{spec intro}

To prove our two main theorems we will need to be able to convert
spectral hypotheses into explicit bounds on geometry.  This section
provides background on the spectrum of the Laplacian for orbifolds, and
establishes facts that will be needed to obtain a diameter bound
in Section \ref{chengchapt}.  Useful references for the material in this
Section are \cite{MR86g:58140} and \cite{ber}.

In this section orbifolds are assumed to be compact and orientable.  The inner 
product on $L^2(O)$ will be indicated
with parentheses,  $(\cdot, \cdot)$.  For vector fields $X$ and
$Y$ on an orbifold $O$, we will use $(X,Y)$ to denote the inner
product $\int_O <X,Y> \ dV$.

Let $O$ be a Riemannian orbifold and let $f$ be a
smooth function on $O$.  The Laplacian $\Delta f$ of $f$ is given by
the Laplacian of lifts of $f$ in the orbifold's local coverings.  More
precisely, lift $f$ to $\tilde{f} = f\circ \pi$ via a coordinate chart $\cc$.  
Let
$g_{ij}$ denote the $\Gamma$-invariant metric on $\tilde{U}$ and $\rho$ =
$\sqrt{\text{det}(g_{ij})}$ as in Section \ref{geod}.  On this local cover
$\Delta\tilde{f}$ is given in the usual way,
\begin{equation*}
\Delta\tilde{f} = \frac{1}{\rho}\sum_{i,j=1}^n
\parttxi(g^{ij}\frac{\partial f}{\partial \tilde{x}^j}\rho).
\end{equation*}

The study of the spectrum of the Laplacian begins with the problem of
finding all of the Laplacian's eigenvalues as it acts on $C^{\lh}(O)$.
That is, we seek all numbers $\lambda$, with multiplicities, that solve
$\Delta f = \lambda f$ for some nontrivial $f \in C^{\lh}(O)$.

The following theorem is proven in \cite{MR94c:58040}.

\begin{theorem} Let $O$ be a Riemannian orbifold.  
\begin{enumerate}\label{ofd spectral theorem}
\item \label{ofd spectral theorem 1} The set of eigenvalues
  $\ld$ in $\Delta f = \lambda f$ consists of an infinite sequence
  $0 \le \overline{\ld_1} < \overline{\ld2} < \overline{\ld3} < \dots \uparrow 
\lh$. 

\item \label{ofd spectral theorem 2} Each eigenvalue $\ld_i$ has finite 
multiplicity.  (Eigenvalues
  will henceforth be listed as $0 \le \ld_1 \le \ld_2 \le \ld_3 \dots
  \uparrow \lh$ with each eigenvalue repeated according to its multiplicity.)
\item \label{ofd spectral theorem 3} There exists an orthonormal basis of $L^2
(O)$ composed of smooth
  eigenfunctions $\vp_1, \vp_2, \varphi_3 \dots$ where
  $\Delta \varphi_i = \ld_i \varphi_i$.
\end{enumerate}
\end{theorem}
The sequence $0 \le \ld_1 \le \ld_2 \le \ld_3 \dots
  \uparrow \lh$ in Theorem \ref{ofd spectral theorem}(\ref{ofd spectral theorem 
2}) is called the spectrum
  of the Laplacian on $O$.  It will be denoted by $Spec(O)$.

The first Sobolev space of a Riemannian orbifold $O$ is obtained
by completing $C^{\infty}(O)$ with respect to the norm associated to
the following inner product,
\begin{equation*}
(f, h)_1 = (f, h) + (\nabla f, \nabla h).
\end{equation*}

\noindent We'll denote the first Sobolev space by $\ch (O)$, and the
associated norm by $||\cdot||_1$.  Note that,
\begin{align*}
C^{\infty}(O) \subset \mathcal{H}(O) \subset L^2(O).
\end{align*}  
Non-smooth elements $u$
of $\ch (O)$ possess first derivatives in the distributional sense.
In analogy with the gradient of a smooth function, these weak
derivatives will be denoted by $\nabla u$. See \cite{MR2001k:58060} for 
information about general orbifold
Sobolev spaces.

A useful tool in spectral geometry is the Rayleigh quotient.  It is defined as 
follows.
\begin{defn}\label{rrq} For $h \in \ch (O)$ with $\int_O h^2 dV \ne 0$ the
  Rayleigh quotient of $h$ is defined by, 
\begin{equation*}
R(h) = \frac{\int_O <\nabla h, \nabla h> dV}{\int_O h^2 dV}.
\end{equation*}
\end{defn}

The proof of Rayleigh's Theorem for the closed eigenvalue problem
extends from the manifold category to the orbifold category without difficulty.

\begin{lemma}(Rayleigh's Theorem for Orbifolds)\label{Rayleigh}  Let
  $O$ be a Riemannian orbifold with
  eigenvalue spectrum $0 \le \ld_1 \le \ld_2 \le \ld_3 \dots \uparrow \lh$.
\begin{enumerate}
\item Then for any $h \in \ch (O)$, with $h \ne 0$, we have $R(h)
  \ge \ld_1$ with equality if and only if $h$ is an eigenfunction of $\ld_1$.
\item Suppose $\{\varphi_1, \varphi_2, \dots \}$ is a complete orthonormal
  basis of $L^2(O)$ with $\varphi_i$ an eigenfunction of $\ld_i$,
  $i=1, 2, 3, \dots$.  If $h \in \ch(O)$ with $h \ne 0$ satisfies $(h,
  \varphi_1) = (h, \varphi_2) = \dots = (h, \varphi_{k-1}) = 0$, then
  $R(h) \ge \ld_k$ with equality if and only if $h$ is an
  eigenfunction of $\lambda_k$.
\end{enumerate}
\end{lemma}

In \cite{MR2001k:58060} it is shown that Weil's asymptotic formula
extends to the orbifold category as well.

\begin{theorem}\label{weil}(Weil's asymptotic formula)  Let $O$ be a
  Riemannian orbifold with eigenvalue spectrum $0 \le \ld_1 \le \ld_2 \le \ld_3
  \dots \uparrow \lh$.  Then for the function $N(\ld) = \sum_{\ld_j
  \le \ld} 1$ we have
\begin{equation*}
N(\ld) \sim (\vol B^n_{0}(1)) (\vol O) \frac{\ld^{n/2}}{(2\pi)^n}
\end{equation*}
as $\ld \uparrow +\lh$.  Here $B^n_{0}(1)$ denotes the $n$-dimensional
unit ball in Euclidean space.
\end{theorem}

\noindent Thus, as with the manifold case, the Laplace spectrum determines an 
orbifold's dimension and
volume.

\section{Obtaining the Diameter Bound}\label{chengchapt}

By applying volume comparison tools in the spectral setting, we derive an upper
diameter bound for an orbifold that relies on spectral information and
the presence of a lower Ricci curvature bound.  With the diameter bound 
established,
an application of the Orbifold Relative Volume Comparison Theorem (\propref
{ofdrvct}) proves the first main theorem.

As in the preceding section, we assume that all orbifolds are
compact and orientable.  Also, we will let $R(\cdot)$ denote
the Rayleigh quotient from Section \ref{spec intro}, Definition \ref{rrq}.

For any open set $U$ in $O$, let $\ch_0(U)$ denote the completion of
$C^{\lh}_0(U)$ in $\ch(U)$.
\begin{defn}
Let $U$ be an arbitrary open set in a Riemannian orbifold $O$.  The
fundamental tone of $U$, denoted $\lambda^*(U)$ is defined by,
\begin{equation*}
\lambda^*(U) = \inf\{R(f)\ : \; f \in \mathcal{H}_0(U), f \ne 0 \}.
\end{equation*}
\end{defn}

The following fact about the fundamental tone will be used in the
proof of the orbifold version of Cheng's Theorem.  Its proof is
identical to that of the manifold version.

\begin{lemma}\label{fund tone fact}
Let $\{U_{\alpha}\}_{\al \in I}$ be a set of domains in a Riemannian
orbifold $O$.  Set $U = \bigcup_{\alpha}
U_{\alpha}$.  Then $\lambda^*(U) \le \inf_{\alpha} \lambda^*(U_{\alpha})$.
\end{lemma}

In what follows let $\scc$ denote the $n$-dimensional simply connected
space form of constant curvature $\kappa$.  Let $B^n_{\kappa}(r)$
denote the ball of radius $r$ in $\scc$, and let $\lambda^n_{\kappa}(r)$
denote the lowest Dirichlet eigenvalue of $B^n_{\kappa}(r)$.

\begin{prop}\label{chth} (Cheng's Theorem for orbifolds.)  Let $O$ be an $n$-dimensional Riemannian orbifold with Ricci curvature bounded
  below by $\kappa(n-1)$, $\kappa$ real.  Then for any $r > 0$ and $p \in
  O$ we have,
\begin{equation*}
\lambda^*(B(p, r)) \le \lambda^n_{\kappa}(r).
\end{equation*}
\end{prop} 

\begin{proof}  
If $p$ is a manifold point in $O$, the manifold proof of Cheng's
Theorem carries over to orbifolds (see \cite{MR86g:58140}). 

Now suppose $p$ is an arbitrary point in $O$, and take $\{p_i\}
\subset (O - \Sigma_O)$ such that $\{p_i\} \rightarrow p$.  Consider the
infinite collection of balls $\{B(p_i,r -
d(p_i,p)\}_{i=1}^\infty$.  Note in particular that $\bigcup_{i=1}^\infty B
(p_i,r -
d(p,p_i))$ is equal to $B(p, r)$.  By Lemma \ref{fund tone fact}
we have 
\begin{equation*}
\lambda^*(B(p,r)) \le \inf_i  \lambda^*(B(p_i, r-d(p_i,p))). 
\end{equation*}

\noindent Since the $p_i$'s are manifold points we can invoke the previous case to obtain,
\begin{equation}\label{chend1}
\lambda^*(B(p,r)) \le \inf_i \lambda^*(B(p_i, r-d(p_i,p)))
\le \inf_i \lambda_{\kappa}^n(r - d(p_i,p)).
\end{equation}
Finally by domain monotonicity of eigenvalues we have,
\begin{equation}\label{chend2}
\inf_i \lambda^n_{\kappa}(r - d(p_i,p)) = \lambda^n_{\kappa}(r).
\end{equation}
Combining lines \ref{chend1} and \ref{chend2} concludes the argument.
\end{proof}

We now adapt a method introduced in \cite{MR93f:53034} to the
orbifold setting.  This method uses spectral data about an orbifold, together 
with a
lower Ricci curvature bound, to obtain an upper bound on the
diameter of the orbifold.  Recall that $\lambda^n_{\kappa}(r)$ denotes
the lowest Dirichlet eigenvalue of $B^n_{\kappa}(r)$.

\begin{prop}\label{diameter bound}
Let $O$ be a compact Riemannian orbifold with Ricci curvature bounded below by
$\kappa(n-1)$, $\kp$ real.  Fix arbitrary constant $r$ greater than zero. Then 
the number of disjoint balls of radius $r$ that can be
placed in $O$ is bounded above by a number that depends only on
$\kappa$ and the number of eigenvalues of $O$ less than or equal to $\lambda_
{\kp}^n(r)$.

In particular the diameter of $O$ is bounded above by a number that
depends only on $Spec(O)$, $\kappa$ and $r$.
\end{prop}

\begin{proof}
 As before, write the eigenvalue spectrum of $O$ as $\lambda_1 \le
\lambda_2 \le \lambda_3 \le \dots \uparrow \infty$.  Choose $\ve > 0$ so that 
no eigenvalues of $O$ lie
between $\lambda_{\kp}^n(r)$ and $\lambda_{\kp}^n(r) + \ve$.  Take a
collection of $N(r)$ pairwise disjoint metric $r$-balls
$B(p_1, r)$, $B(p_2, r)$, $\dots$, $B(p_{N(r)},r)$ in $O$.  By Cheng's
Theorem (Proposition \ref{chth}) we have for each $i$ a function
$f^i \in \ch_0(B(p_i,r))$ such that $R(f^i) \le \lambda_{\kp}^n(r) +
\ve$.

Because $\ch_0(B(p_i,r))$ is the closure of
$C^{\infty}_0(B(p_i,r))$ with respect to $||\cdot||_1$ we can
find for each $i$ a sequence $\{h_j^i\}_{j=1}^{\lh} \subset C^{\infty}_0(B
(p_i,r))$
that converges to $f^i$.  By the continuity of $R: \ch(O) \rightarrow
\mathbf{R}$ we know additionally that $R(h^i_j) \rightarrow R(f^i)$
as $j \rightarrow \infty$.  In particular for $\ve'>0$ arbitrary we
can choose integers $k(i)$ large enough that $|R(h^i_{k(i)}) - R(f^i)|
< \ve'$ for each $i$.  

Extend each $h^i_{k(i)}$ to all of $O$ by setting it equal to zero
off of $B(p_i, r)$.  Now $(h^i_{k(i)},h^j_{k(j)}) = 0$ for $i \ne j$ as the 
supports of these functions
are disjoint.  To arrange that the collection $\{h^i_{k(i)}\}_{i=1}^\infty$
is orthonormal replace each $h^i_{k(i)}$ with $\overline{h}^i_{k(i)}=\frac{h^i_
{k(i)}}{|h^i_{k(i)}|}$.

Pick  $\vp_1, \vp_2, \dots, \vp_{N(r) - 1} \in
L^2(O)$ which are orthonormal and which are eigenfunctions for
$\lambda_1, \lambda_2, \dots , \lambda_{N(r)-1}$ respectively.
There exist $\al_1, \al_2, \dots, \al_{N(r)}$, not all
zero, such that,
\begin{align*}
\Sigma^{N(r)}_{l=1} \al_l(\overline{h}^l_{k(l)}, \varphi_m) = 0,
\end{align*} 
for $m = 1, 2, \dots , N(r)-1$.  
Setting $\psi = \Sigma^{N(r)}_{l=1} \al_l \overline{h}_{k(l)}^l$, Rayleigh's 
Theorem yields, 
\begin{align*}
\lambda_{N(r)}|\psi|^2 &\le (\nabla \psi, \nabla \psi) \\
&= (\Sigma^{N(r)}_{l=1} \al_l \nabla \overline{h}_{k(l)}^l, \Sigma^{N(r)}_{l=1} 
\al_l \nabla \overline{h}_{k(l)}^l)\\
&= \Sigma^{N(r)}_{l,s=1} \al_l \al_s (\nabla \overline{h}_{k(l)}^l, \nabla 
\overline{h}_{k(s)}^s) \\
&= \Sigma^{N(r)}_{l=1} \al_l^2 (\nabla \overline{h}_{k(l)}^l, \nabla \overline
{h}_{k(l)}^l) \\
&= \Sigma^{N(r)}_{l=1} \al_l^2 \frac{(\nabla h_{k(l)}^l, \nabla h_{k(l)}^l)}{|h_
{k(l)}^l|^2}\\
&= \Sigma^{N(r)}_{l=1} \al_l^2 R(h^l_{k(l)}). 
\end{align*}
Now observe that $|\psi|^2 = \Sigma^{N(r)}_{l=1} \al_l^2$.  The
calculation above then implies,
\begin{equation*}
\lambda_{N(r)} \Sigma^{N(r)}_{l=1} \al_l^2 \le \Sigma^{N(r)}_{l=1} \al_l^2 R
(h^l_{k(l)}).
\end{equation*}
By our choice of $k(l)$ we have,
\begin{equation*}
\lambda_{N(r)} \Sigma^{N(r)}_{l=1} \al_l^2 \le \Sigma^{N(r)}_{l=1} \al_l^2 (R
(f^l)+\ve').
\end{equation*}
And our choice of $f^l$ gives,
\begin{equation*}
\lambda_{N(r)} \Sigma^{N(r)}_{l=1} \al_l^2 \le \Sigma^{N(r)}_{l=1} \al_l^2 
(\lambda^n_{\kappa}(r) + \ve + \ve').
\end{equation*}
Since we know at least one $\al_l$ is nonzero, we can divide both
sides by $\Sigma^{N(r)}_{l=1} \al_l^2$ to obtain $\lambda_{N(r)}
\le \lambda^n_{\kappa}(r) + \ve + \ve'$.  Letting $\ve'$ go to zero
simplifies the right hand side further and we have,
\begin{equation*}
\lambda_{N(r)} \le \lambda^n_{\kappa}(r) + \ve.
\end{equation*}
Because $\ve$ was chosen so that no eigenvalues of $O$ appeared in
$(\lambda^n_{\kappa}(r), \lambda^n_{\kappa}(r) + \ve)$, we can
conclude that $\lambda_{N(r)} \le \lambda^n_{\kappa}(r)$.

We now obtain the diameter bound.  Let $\rho$ be the largest number
so that $\lambda_{\rho}(O) \le \lambda^n_{\kappa}(r)$.    Then
$\lambda_{N(r)} \le \lambda^n_{\kappa}(r)$ implies that $N(r) \le \rho$.  Thus 
the number of disjoint $r$-balls is bounded by the spectral
invariant $\rho$.  Now an orbifold of diameter $d$ contains at
least $[d/2r]$ disjoint $r$-balls, where $[\cdot]$ denotes the greatest integer
function.  Thus $d$ must satisfy $[d/2r] \le \rho$.  This gives an
upper bound on the diameter of $O$ which depends only on $r$,
$\kappa$ and $Spec(O)$.
\end{proof}

We are now prepared to prove our first main result.  

\medskip

\noindent \bf{Main Theorem 1:} \it  
Let $\cs$ be a collection of isospectral Riemannian
orbifolds that share a uniform lower bound $\kp(n-1)$, $\kp$ real, on Ricci 
curvature.  Then there are only finitely
many possible isotropy types, up to isomorphism, for points in an
orbifold in $\cs$. \rm

\begin{proof}
It is shown above that isospectral families of
orbifolds which share a uniform lower Ricci curvature bound also
share an upper diameter bound.  Let $D>0$ be the upper bound for the
diameter of orbifolds in $\cs$.  By Weil's asymptotic formula, the 
isospectrality
of the orbifolds in $\cs$ implies that they all have the same
dimension $n$ and the same volume $v>0$.  

Let $O$ be an orbifold in $\cs$ and take $p \in O$.  As before let
$B^n_{\kp}(r)$ denote the $r$-ball in the simply connected,
$n$-dimensional space form of constant curvature $\kp$.  For $R>r\ge 0$
we have by Proposition \ref{ofdrvct},
\begin{equation*}
\frac{\vol B(p,r)}{B^n_{\kp}(r)} \ge \frac{\vol B(p,R)}{B^n_{\kp}(R)}.
\end{equation*}
Letting $R = D$ in this inequality gives,
\begin{equation*}
\frac{\vol B(p,r)}{B^n_{\kp}(r)} \ge \frac{\vol O}{B^n_{\kp}(D)} = \frac{v}{B^n_
{\kp}(D)}.
\end{equation*}
Again applying Proposition \ref{ofdrvct} we take the limit as $r
\rightarrow 0$ to obtain,
\begin{equation*}
\frac{1}{|\Gamma_p|} = \lim_{r \rightarrow 0}\frac{\vol B(p,r)}{B^n_{\kp}(r)} 
\ge
\frac{v}{B^n_{\kp}(D)}.
\end{equation*}
We conclude for any point in any orbifold in $\cs$, the isotropy group
of that point has order less than or equal to the universal constant
$B^n_{\kp}(D)/v$.  This implies that the isotropy group of the
point can have one of only finitely many possible isomorphism types. 
\end{proof}

Consider the collection of all closed, connected Riemannian
$n$-orbifolds with a lower bound $\kappa(n-1)$, $\kp$ real, on Ricci
curvature, a lower bound $v>0$ on volume, and an upper bound $D>0$ on 
diameter.  A similar argument to the
one above shows that there are only finitely many possible isotropy
types for points in an orbifold in this collection.

\section{Spectral Bounds on Isolated Singular Points}\label{main lemma chap}

This section begins by extending a technical result from \cite{MR90a:53044} to 
the
orbifold setting.  Assume the orbifolds under consideration are compact and 
orientable.  As in Section \ref{ofdcompgeom}, we will use $B^{\fa}(p,r)$ to 
denote
the cone of radius $r$ at point $p$ in an orbifold with directions
given by $\fa \in S_pO$.  Following \cite{MR90a:53044} we will use the symbol
$\oct$ to denote the collection of all closed, connected
$n$-dimensional Riemannian orbifolds with volume bounded below by $v>0$, 
sectional
curvature bounded below by $\kp \in \mathbf{R}$, and with diameter bounded 
above by $D>0$.  The
subcollection of orbifolds in $\oct$ with only isolated singularities
will be denoted by $isol\oct$.

Suppose $O$ is a complete orbifold and $K$ is a compact subset of
$O$.  Let $\apK$ denote
the set of unit tangent vectors at $p$ which are the velocity vectors
of segments running from $p$ to $K$.  The set $\apK$ is called the
\emph{set of directions from $p$ to $K$}.  

For subset $\fa$ of the unit $n$-sphere, $S^n$, we write,
\begin{align*}
\fa(\theta) &= \{v \in S^n : \angle(\fa,v) < \theta\} \\
\fa'(\theta) &= \{v \in S^n : \angle(\fa,v) \ge \theta\}.
\end{align*}

\begin{lemma} \label{appdxnolemma} Suppose for some $\alpha \in [0,
  \frac{\pi}{2}]$, a finite subset $A$ of $S^n$ satisfies,
\begin{equation*}
\overline{A\bigl(\pt + \al\bigr)} = S^n.
\end{equation*}
Let $\tilde{\mathfrak{a}}_{\al} \subset S^n$ consist of two vectors situated
  at an angle of $\pi - 2\al$ from each other.  Then, using the standard
  volume on $S^n$, we have:
\begin{equation*}
\text{Vol } A(\theta) \ge
  \text{Vol } \tilde{\mathfrak{a}}_{\al}(\theta)
\end{equation*}
for all $\theta$ greater than
  or equal to zero.
\end{lemma}

\begin{proof}
See the appendix in \cite{MR90a:53044}.
\end{proof}

\begin{lemma}\label{main lemma 1} Let $O \in isol\oct$ and $p, q
  \in O$.   There exist $\alpha \in (0,
  \frac{\pi}{2})$ and $r > 0$ such that if,
\begin{equation*}
\mathfrak{d}_{pq}(\frac{\pi}{2} + \alpha) =
S_pO, \ \text{and} \ \mathfrak{d}_{qp}(\frac{\pi}{2} + \alpha) = S_qO,
\end{equation*}
then $d(p, q) \ge r$.  The constants $\al$ and $r$ depend only on $\kp$,
$D$, $v$ and $n$.
\end{lemma}

A remark on the positive curvature case is necessary before proving
this lemma.  If $\kp > 0$ then the Bonnet-Myers Theorem implies that for $O \in
\oct$, the manifold $O
- \Sigma_O$ has diameter less than or equal to $\pi/\sqrt{\kp}$.  Thus $O$
itself satisfies this diameter bound.  So in the positive curvature
case we can assume $D \le \pi/\sqrt{\kp}$.  In particular orbifolds in $\oct$
satisfy the hypotheses of the Volume Comparison Theorem for cones in
orbifolds (Proposition \ref{ofdcvct}), which will be used below.

\begin{proof}(Lemma \ref{main lemma 1})
For a parameter $\alpha \in (0, \frac{\pi}{2})$ let
$\tilde{\mathfrak{a}}$ be a subset of $S^{n-1}$ consisting of two vectors,
$v$ and $w$, for which $\angle(v, w) = \pi -2\alpha$.  Let $\bar{p}$ be an 
element of $\scc$, the simply connected complete $n$-dimensional space form of 
constant curvature $\kp$.  We specify $\alpha$ by choosing it as an
element of $(0, \frac{\pi}{2})$ such that:
\begin{equation*}
\vol B^{\tilde{\mathfrak{a}}'(\frac{\pi}{2} - \alpha)}(\bar{p},D)< \frac{v}{6}.
\end{equation*}

Suppose we have points $p$ and $q$ in $O$ for which,
\begin{equation*}
\mathfrak{d}_{pq}(\frac{\pi}{2} + \alpha) =
S_pO, \ \text{and} \ \ \mathfrak{d}_{qp}(\frac{\pi}{2} + \alpha) = S_qO.
\end{equation*}
\begin{figure}
\scalebox{.5}{\includegraphics{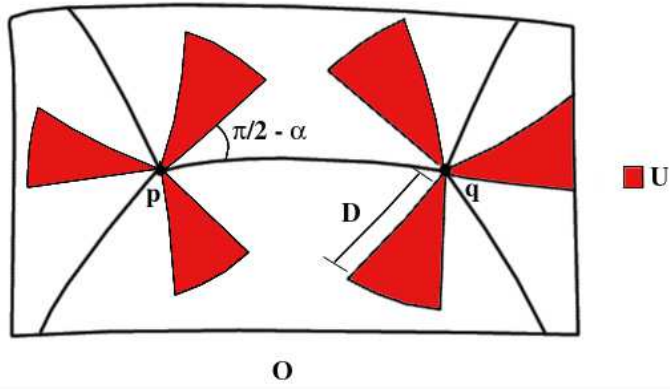}}
\caption{The Set $U$}\label{ufig}
\end{figure}

\noindent Now, $S_pO$ and $S_qO$ are compact so we can take finite
  subsets $\mathfrak{d}_p$ and $\mathfrak{d}_q$ of ${\mathfrak{d}}_{pq}$ and
  ${\mathfrak{d}}_{qp}$, respectively, so that 
\begin{equation*}
\mathfrak{d}_p(\frac{\pi}{2} + \alpha) = S_pO, \ \text{and} \ \ \mathfrak{d}_q
(\frac{\pi}{2} + \alpha) = S_qO 
\end{equation*}
as well.
Lifting these sets gives,
\begin{equation*}
\tilde{\mathfrak{d}}_p(\frac{\pi}{2} + \alpha) =
S_{\tilde{p}}\widetilde{U}_p, \ \text{and} \ \ \tilde{\mathfrak{d}}_q(\frac{\pi}
{2} + \alpha) = S_{\tilde{q}}\widetilde{U}_q.
\end{equation*} 

\noindent Because of this we can use Lemma \ref{appdxnolemma}
to conclude that,
\begin{equation*}
\vol(\tilde{\mathfrak{d}}_p'(\frac{\pi}{2} - \alpha)) \le
\vol(\tilde{\mathfrak{a}}'(\frac{\pi}{2} - \alpha)), \ \text{and} \ \ \vol
(\tilde{\mathfrak{d}}_q'(\frac{\pi}{2} - \alpha)) \le
\vol(\tilde{\mathfrak{a}}'(\frac{\pi}{2} - \alpha)).
\end{equation*}

Let $U$ be the subset of $O$ given by,
\begin{align*}
U &= B^{\inr(\mathfrak{d}_p'(\frac{\pi}{2} -
\alpha))}(p,D) \cup B^{\inr(\mathfrak{d}_q'(\frac{\pi}{2} -
\alpha))}(q,D),
\end{align*}
where $\inr(\mathfrak{d}_p'(\frac{\pi}{2} -
\alpha))$ denotes the interior of $\mathfrak{d}_p'(\frac{\pi}{2} -
\alpha)$, and $\inr(\mathfrak{d}_q'(\frac{\pi}{2} -
\alpha))$ denotes the interior of $\mathfrak{d}_q'(\frac{\pi}{2} -
\alpha)$.  A sketch of the set $U$ is given in \figref{ufig}.  The lines
emanating from $p$ and $q$ indicate the segments between these
points.  The shaded regions are the cones that form $U$.

Let $I:S_{\tp}\wtu_p \rightarrow S_{\op}\scc$ and $J:S_{\tq}\widetilde{U}_q
\rightarrow S_{\op}\scc$ be linear isometries.  Then using Proposition \ref
{ofdcvct} we have that $\vol(U) < v/3$, as:
\begin{align*}
\vol U &\le \vol B^{\inr(\mathfrak{d}_p'(\frac{\pi}{2} -
\alpha))}(p,D) + \vol B^{\inr(\mathfrak{d}_q'(\frac{\pi}{2} -
\alpha))}(q,D) \\
&\le \vol B^{I(\inr(\tilde{\mathfrak{d}}_p'(\frac{\pi}{2} -
\alpha)))}(\op,D) + \vol B^{J(\inr(\tilde{\mathfrak{d}}_q'(\frac{\pi}{2} -
\alpha)))}(\oq,D) \\
&= \vol \e_{\bar{p}}[0, D]I(\inr(\tilde{\mathfrak{d}}_p'(\frac{\pi}{2} -
\alpha))) + \vol \e_{\bar{q}}[0,
D]J(\inr(\tilde{\mathfrak{d}}_q'(\frac{\pi}{2} - \alpha)))\\
&\le \vol \e_{\bar{p}}[0, D]I(\tilde{\mathfrak{d}}_p'(\frac{\pi}{2} -
\alpha)) + \vol \e_{\bar{q}}[0,
D]J(\tilde{\mathfrak{d}}_q'(\frac{\pi}{2} - \alpha))\\
&\le 2 \vol \e_{\bar{p}}[0, D]\tilde{\mathfrak{a}}'(\frac{\pi}{2} -
\alpha) \\
&< \frac{v}{3}. 
\end{align*}

We are ready to specify the constant $r$ required by the statement of
the Lemma.  First choose $l > 0$ so that $\vol \ccb(l) <
v/3$ in $M^n_k$.  Note that the Orbifold Relative Volume
Comparison Theorem (Proposition \ref{ofdrvct}) implies,
\begin{equation*}
\vol(B(p, l)) < \frac{v}{3},\ \text{and} \ \ \vol(B(q, l)) < \frac{v}{3}.
\end{equation*}
Let $(c_1; c_2; c_3)$ denote a geodesic triangle in $\scc$ with sides
$c_1$, $c_2$ and $c_3$.  In triangle $(c_1; c_2; c_3)$ the angle opposite side
$c_i$ will be denoted by $\theta_i$.  For the $\alpha \in (0,
\frac{\pi}{2})$ and $l > 0$ determined above, there exists an $r > 0$ such that 
for all geodesic triangles $(c_1;
c_2; c_3)$ in $\scc$ with $\theta_1 \le \frac{\pi}{2} - \alpha$,
$L(c_3) \ge l$ and $L(c_2) < r < l$, we also have $L(c_1) < L(c_3)$.
\begin{figure}
\scalebox{.5}{\includegraphics{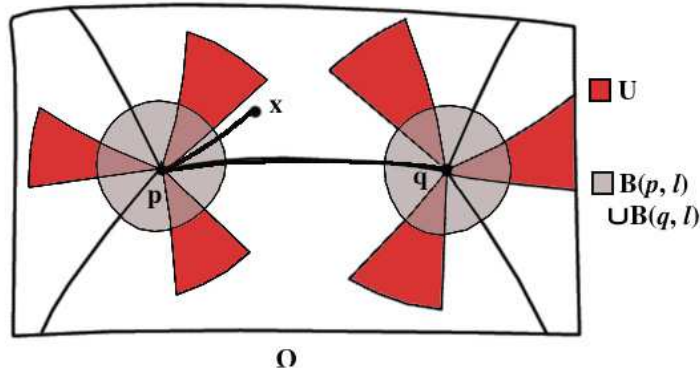}}
\caption{Hinge in $O$}\label{hinge}
\end{figure}
\begin{figure}
\scalebox{.5}{\includegraphics{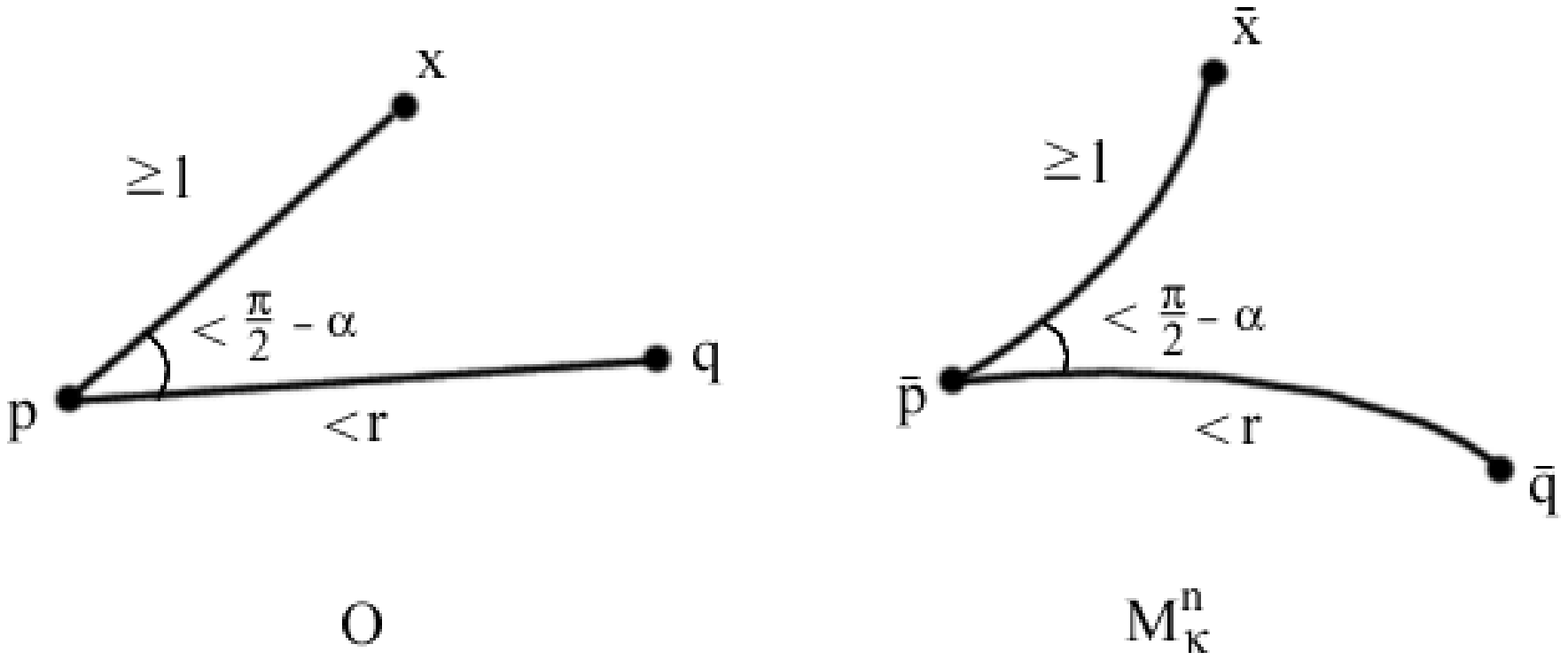}}
\caption{Hinge in $O$; Comparison Hinge in $\scc$}\label{contrahinge}
\end{figure}

We finish the proof by nested contradiction arguments.  That is, we
will show that if $p$ and $q$ satisfy the hypotheses of the Lemma, and
$d(p, q) < r$, then the sets $U$, $B(p, l)$ and $B(q, l)$ cover $O$.  However if
  these sets cover $O$ we have, 
\begin{equation*}
v \le \vol(O) \le \vol(U) + \vol(B(p, l)) + \vol(B(q,
l)) < v.
\end{equation*}
Since this is a contradiction, once we show that $U$, $B(p,
l)$, and $B(q, l)$ cover $O$ we can conclude that $d(p, q) \ge r$.

To show that $U$, $B(p, l)$, and $B(q, l)$ cover $O$ we argue again by
contradiction.  Suppose they fail to cover and we can find a point $x$ in
$O - (U \cup B(p, l) \cup B(q, l))$.  Set up a hinge with angle
at $p$ terminating at $x$ and $q$ so that the leg of the hinge from
$p$ to $q$ is a segment, and so that the hinge angle is less than
$\frac{\pi}{2} - \alpha$.  \figref{hinge} gives a sketch of
this hinge in $O$.
Using Toponogov's Theorem for hinges in orbifolds, form a
comparison hinge in $\scc$ with angle at $\bar{p}$ terminating at
$\bar{x}$ and $\bar{q}$.  \figref{contrahinge} illustrates both the
original hinge in $O$ and the comparison hinge in $\scc$. 

By Toponogov's Theorem we know $d(x, q) \le d(\bar{x},\bar{q})$.  In addition 
our
setup implies that $d(\bar{p},\bar{q}) < r$, $d(\op, \bar{x}) \ge l$, and
the angle formed by the comparison hinge is less than $\frac{\pi}{2} -
\alpha$.  So by our choice of $r$ we have $d(\bar{x},\bar{q}) <
d(\bar{x},\bar{p})$.  Putting these observations together yields, 
\begin{equation*}
d(x, q) \le d(\bar{x}, \bar{q}) < d(\bar{x}, \bar{p}) = d(x, p),
\end{equation*}
thus $d(x, q) < d(x, p)$.  A similar argument based at $q$ yields the
contradictory statement $d(x, p) < d(x, q)$, completing the proof.
\end{proof}

We now use \lemref{main lemma 1} to bound the number
of singular points that can appear in an orbifold in $isol\oct$.  This
in turn will lead to our second main theorem.

Fix $\varepsilon > 0$.  A \emph{minimal $\varepsilon$-net}
 in a compact, connected metric space $X$ is an ordered set of points
 $p_1, p_2, \dots, p_N$ with the following two properties.  First, the
 open balls $B(p_i, \varepsilon)$, $i = 1, 2, \dots, N$, cover $X$.
 Second, the open balls $B(p_i, \varepsilon/2)$, $i = 1, 2, \dots, N$
 are disjoint.  The fact that for any $\ve > 0$ we can find a minimal $\ve$-net 
in $X$
 is well known.
 
\begin{prop}\label{cpbound}  There is a positive integer $C(D,v,\kp,n)$ for 
which no orbifold $O$ in the family $isol\oct$ has more than $C$ singular 
points.
\end{prop}

\begin{proof}
Suppose $O \in isol \oct$, and let $\al$ and $r$ be as in Lemma \ref{main lemma 
1}.  Take $p \in \Sigma_O$ and let $(U, \wtu/\Ga_p, \pi)$ be a fundamental
coordinate chart about $p$.  Also, let $\tp$ denote the point in $\wtu$
which projects to $p$ under $\pi$.  The set of lifts of a vector $v
\in S_pO$ is the orbit $\Gamma_{p*}\tv$ of any vector $\tv \in
S_{\tp}\wtu$ for which $\pi_{*_{\tp}} \tv = v$.  We will first show
that $\Gamma_{p*}\tv$ does not lie in any open hemisphere of
$S_{\tp}\wtu$.  With this established we can then appeal to Lemma
\ref{main lemma 1} to conclude that the distance between two singular
points in $O$ will always be greater than $r$.  This in turn will be
used to obtain the universal upper bound on the number of singular
points in $O$.

Because $p$ is an isolated singularity, elements of $\Gamma_{p*}$ act on $S_
{\tp}\wtu$ without fixed points.  Thus the possible quotients $S_{\tp}\wtu / 
\Gamma_{p*}$ are actually all spherical space forms.  Spherical space forms 
obtained as quotients of the sphere by finite groups of orthogonal 
transformations are well understood.  See \cite{MR49:7958} for example.  In 
even dimensions the only non-trivial quotient is projective space, obtained as 
the quotient of $S^{2m}$ by the antipodal map.  Since the orbits under the 
antipodal map consist of pairs of antipodal points, its clear that no orbit is 
contained in an open hemisphere.  

Odd-dimensional spherical space forms, however, can arise in many
ways.  In this situation it will suffice to consider only those that
are quotients of an odd dimensional sphere by the action of a cyclic
group.  This is because if we take an element $\gamma_{*_{\tp}} \in \Gamma_{p*}
$ of order $l$,
to show $\Gamma_{p*}\tv$ is not contained in an open hemisphere it
suffices to show that $\{\tv, \gamma_{*_{\tp}} \tv, \gamma_{*_{\tp}}^2 \tv, 
\dots,
\gamma_{*_{\tp}}^{l-1} \tv\}$ is not contained in any open hemisphere.  

Suppose $\Gamma \le O(2m)$ is cyclic and generated by $\gamma \in
\Gamma$ of order $l$.  Viewing $\mathbf{R}^{2m}$ as $\mathbf{C}^{m}$,
element $\gamma$ 
can be expressed as:

\medskip

\begin{center}

$\begin{bmatrix}
e^{2\pi i/l} & 0 & 0 & \dots & 0 \\
0 & e^{2\pi i a_1 /l} & 0 & \dots & 0 \\
0 & 0 & e^{2\pi i a_2 /l} & \dots & 0 \\
\vdots & \vdots & \vdots & \ddots & \vdots \\
0 & 0 & 0 & \dots & e^{2\pi i a_{m-1}/l} \\ 
\end{bmatrix}$

\end{center}

\medskip

\noindent for $a_1, a_2, \dots, a_{m-1} \in \mathbf{R}$ each relatively prime 
to $l$.  Thus the orbit of a vector $z = (z_1, z_2, \dots, z_m) \in S^{2m-1}$ 
has the following form:
\begin{align*}
\{(e^{2\pi i/l}z_1, e^{2\pi i a_1 /l}z_2, \dots, e^{2\pi i a_{m-1}/l}z_m), &(e^
{2\cdot2\pi i/l}z_1, e^{2\cdot2\pi i a_1 /l}z_2, \dots, e^{2\cdot2\pi i a_{m-
1}/l}z_m), \dots, \\
&(e^{2(l-1)\pi i/l}z_1, e^{2(l-1)\pi i a_1 /l}z_2, \dots, e^{2(l-1)\pi i a_{m-
1}/l}z_m)\}.
\end{align*}
If we sum together all of the orbits of $z$ under $\gamma$ we get the following 
vector in $\mathbf{R}^{2m}$:
\begin{equation}\label{zerovector}
(\Sigma^{l-1}_{k=0}e^{2\pi k i/l}z_1, \Sigma^{l-1}_{k=0}e^{2\pi i a_1k /l}z_2, 
\dots, \Sigma^{l-1}_{k=0}e^{2\pi i a_{m-1}k /l }z_m).
\end{equation}
By showing that this vector is actually the zero vector we will be
able to conclude that  $\{z, \gamma z, \gamma^2 z, \dots, \gamma^{l-1}
z\}$ does not lie in any open hemisphere.  

To see that the vector in line \ref{zerovector} is the zero vector consider the 
$s^{th}$ entry, 
\begin{equation*}
\Sigma^{l-1}_{k=0}e^{2\pi i a_{s-1}k /l}z_s.
\end{equation*}
Since $a_{s-1}$ and $l$ are relatively prime, the set $\{e^{2\pi i
  a_{s-1}k /l}\}^{l-1}_{k=0}$ consists of $l^{\rm th}$ roots of unity.  Because 
the sum of the $l^{\rm
  th}$ roots of unity is zero, we can conclude that this entry
vanishes.  

Now consider points $p$ and $q$ in the singular set of $O$.  Because
$O$ is complete we know that $p$ and $q$ are joined by at least one
segment.  Thus the set of directions from $p$ to $q$ contains at least
one vector, namely the initial vector $v$ of the segment from $p$ to
$q$.  Moreover $\apq(\pt + \al) = S_pO$.  For if this were not the
case we could find $w \in S_pO$ with $\angle(v,w) \ge \pt + \al$.
However, this implies that if we let $\tilde{w}$ be a fixed lift of
$w$ in the covering sphere $S_{\tp}\wtu$, then the orbit of a lift of
$v$ is going to remain within the open hemisphere about $-\tilde{w}$.
This contradicts our conclusions above.  A similar argument shows that
$\aqp(\pt + \al) = S_qO$.  Thus by Lemma \ref{main lemma 1} we know that $d(p, 
q) \ge r$.  

The proof concludes with a volume comparison argument.  Let $\{x_1, x_2, \dots, 
x_N\}$ be a minimal $(r/2)$-net in $O$.  Recall that for $p \in O$ and $S \ge s 
\ge 0$, Proposition \ref{ofdrvct} gives,
\begin{equation}\label{volfact1}
\frac{\vol B^n_{\kp}(s)}{\vol B^n_{\kp}(S)} \le \frac{\vol B(p,s)}{\vol B(p,S)}.
\end{equation}

Without loss of generality suppose that $B(x_1, r/4)$ is the minimal volume 
$(r/4)$-ball in our net.  Then using the fact that the $(r/4)$-balls are 
disjoint we have,
\begin{equation*}
N\vol B(x_1, r/4) \le \Sigma_{i=1}^{N}\vol B(x_i, r/4) \le \vol O.
\end{equation*}
Thus $\vol B(x_1, r/4) \le \vol O/N$. 

Now apply line \ref{volfact1} to balls about $x_1$ with $s = r/4$ and $S = D$.  
This yields,
\begin{equation}\label{volfact2}
\frac{\vol B^n_{\kp}(r/4)}{\vol B^n_{\kp}(D)} \le \frac{\vol B(x_1,r/4)}{\vol B
(x_1,D)}.
\end{equation}
Using $\vol B(x_1,D) = \vol O$ and $\vol B(x_1, r/4) \le \vol O/N$ we find that 
line \ref{volfact2} becomes,
\begin{equation*}
\frac{\vol B^n_{\kp}(r/4)}{\vol B^n_{\kp}(D)} \le \frac{1}{N}.
\end{equation*}
Thus we see that the number of elements in our minimal $(r/2)$-net is
bounded above by the universal constant $\vol B^n_{\kp}(D)/\vol B^n_{\kp}(r/4)$.

The singular points are all at least $r$-apart from each other, so there can be 
at most one singular point per $(r/2)$-ball in our net.  Thus the bound on the 
number of elements in our net is also a bound on the number of singular points 
in $O$.
\end{proof}

Our second main result is a corollary to this proposition.  

\medskip

\noindent \bf{Main Theorem 2:} \it  
Let $isol \cs$ be a collection of isospectral Riemannian orbifolds with only
isolated singularities that share a uniform lower bound $\kp \in
\mathbf{R}$ on sectional curvature.  Then there is an upper
bound on the number of singular points in any orbifold, $O$, in $isol  \cs$
depending only on $Spec(O)$ and $\kp$. \rm

\begin{proof}
The argument begins in the same manner as that in the proof of Main
theorem 1.  Because these orbifolds are isospectral, and satisfy a
lower bound on sectional curvature, we can conclude that they also
share an upper diameter bound.  By
Weil's asymptotic formula, we know that all orbifolds in $isol \cs$ have the 
same
volume and dimension.  Therefore the family $isol \cs$ satisfies the
hypotheses of Proposition \ref{cpbound} and the theorem follows.
\end{proof}


\bibliography{bibfile}

\begin{thebibliography}{BPP92}

\bibitem[B{\'e}r86]{ber}
Pierre~H. B{\'e}rard.
\newblock {\em Spectral geometry: direct and inverse problems}.
\newblock Springer-Verlag, Berlin, 1986.
\newblock With appendixes by G\'erard Besson, and by B\'erard and Marcel
  Berger.

\bibitem[Bor93]{MR94d:53053}
Joseph~E. Borzellino.
\newblock Orbifolds of maximal diameter.
\newblock {\em Indiana Univ. Math. J.}, 42(1):37--53, 1993.

\bibitem[BPP92]{MR93f:53034}
Robert Brooks, Peter Perry, and Peter Petersen, V.
\newblock Compactness and finiteness theorems for isospectral manifolds.
\newblock {\em J. Reine Angew. Math.}, 426:67--89, 1992.

\bibitem[Cha84]{MR86g:58140}
Isaac Chavel.
\newblock {\em Eigenvalues in {R}iemannian geometry}.
\newblock Academic Press Inc., Orlando, FL, 1984.
\newblock Including a chapter by Burton Randol, With an appendix by Jozef
  Dodziuk.

\bibitem[Chi93]{MR94c:58040}
Yuan-Jen Chiang.
\newblock Spectral geometry of ${V}$-manifolds and its application to harmonic
  maps.
\newblock In {\em Differential geometry: partial differential equations on
  manifolds (Los Angeles, CA, 1990)}, pages 93--99. Amer. Math. Soc.,
  Providence, RI, 1993.

\bibitem[Far01]{MR2001k:58060}
Carla Farsi.
\newblock Orbifold spectral theory.
\newblock {\em Rocky Mountain J. Math.}, 31(1):215--235, 2001.

\bibitem[GP88]{MR90a:53044}
Karsten Grove and Peter Petersen, V.
\newblock Bounding homotopy types by geometry.
\newblock {\em Ann. of Math. (2)}, 128(1):195--206, 1988.

\bibitem[Gro99]{MR2000d:53065}
Misha Gromov.
\newblock {\em Metric structures for {R}iemannian and non-{R}iemannian spaces}.
\newblock Birkh\"auser Boston Inc., Boston, MA, 1999.
\newblock Based on the 1981 French original [MR 85e:53051], With appendices by
  M.\ Katz, P.\ Pansu and S.\ Semmes, Translated from the French by Sean
  Michael Bates.

\bibitem[Pet98]{MR98m:53001}
Peter Petersen.
\newblock {\em Riemannian geometry}.
\newblock Springer-Verlag, New York, 1998.

\bibitem[Sat56]{MR18:144a}
I.~Satake.
\newblock On a generalization of the notion of manifold.
\newblock {\em Proc. Nat. Acad. Sci. U.S.A.}, 42:359--363, 1956.

\bibitem[Sat57]{MR20:2022}
Ichir\^o Satake.
\newblock The {G}auss-{B}onnet theorem for ${V}$-manifolds.
\newblock {\em J. Math. Soc. Japan}, 9:464--492, 1957.

\bibitem[Shi93]{MR96c:53064}
Katsuhiro Shiohama.
\newblock {\em An introduction to the geometry of {A}lexandrov spaces}.
\newblock Seoul National University Research Institute of Mathematics Global
  Analysis Research Center, Seoul, 1993.

\bibitem[Thu78]{thur}
William Thurston.
\newblock {\em The Geometry and Topology of 3-Manifolds}.
\newblock Lecture Notes, Princeton University Math. Dept., 1978.

\bibitem[Wol74]{MR49:7958}
Joseph~A. Wolf.
\newblock {\em Spaces of constant curvature}.
\newblock Publish or Perish Inc., Boston, Mass., third edition, 1974.

\end{thebibliography}

\end{document}